\documentclass[twoside,11pt]{article}

\usepackage{amsmath,amsfonts}

\begin{document}
\title{\bf Integral Solutions to Linear
Indeterminate Equation }
\author{C{\sc
hangjiang} Z{\sc hu} \\
The Hubei Key Laboratory of Mathematical Sciences, \\
School of Mathematics and Statistics,\\
Huazhong Normal University, Wuhan 430079, P.R. China.}

\date{}
\maketitle

\vskip 1cm

\noindent {\bf Abstract:} \ \ In this paper,
using Euler's function, we give a formula of all integral
solutions to linear indeterminate equation with $s$-variables
$a_1x_1+a_2x_2+\cdots+a_sx_s=n$.
It is a explicit formula of the coefficients
$a_1$, $a_2$, $\cdots$, $a_s$ and the free term $n$.

\bigbreak \noindent {\bf Key words}: \ \ Linear indeterminate
equation, Euler's function, integral solution.

\bigbreak \noindent{\bf 2000 AMS Subject Classification}: 11D04,
11D72.

\section*{1.   Introduction and Main Theorem}

In this paper, we consider the integral
solutions to linear indeterminate equation with $s$-variables
$$
a_1x_1+a_2x_2+\cdots+a_sx_s=n.
\eqno(1.1)
$$
It is well-known that there exist the integral solutions of (1.1) if
and only if
$$
(a_1,a_2, \cdots, a_s)|n.
\eqno(1.2)
$$
Under the assumption (1.2), if we can obtain a special solution of
(1.1) by applying the mutual division, fraction, parameter methods,
etc., then the all integral solutions to (1.1) can be represented by
using the special solution obtained above and $s-1$ parameters $t_1,
\ t_2, \ \cdots, \ t_{s-1}$. However, the methods seeking above
special solution are too complicated to lose availability in many
problems. For example, it is very difficult to obtain a special
solution to the following simple indeterminate equation with $s=2$
$$
2^mx+3^ny=1,
\eqno(1.3)
$$
where $m$ and $n$ are positive integers. Therefore, it is very
important and interesting to seek a formula of all integral
solutions to (1.1). In this paper, using Euler's function, we give a
formula of all integral solutions to (1.1), which is a explicit
function of the coefficients $a_1$, $a_2$, $\cdots$, $a_s$ and the
free term $n$.

To state our result, let $ (a_1, a_2, \cdots, a_s)=d, $ $n=dn_1$,
$(a_1, a_2)=d_2$, $(d_2, a_3)=d_3, \ \cdots$, $(d_{s-1},
a_s)=d_s=d$, $a_1=d_2\bar{a}_1$, $a_2=d_2\bar{a}_2, \ \cdots$,
$a_s=d_s\bar{a}_s$, $d_2=d_3\bar{d}_2$, $d_3=d_4\bar{d}_3, \
\cdots$, and $d_{s-1}=d_s\bar{d}_{s-1}$.

Then
$$
(\bar{a}_1, \bar{a}_2)=1, \ \ (\bar{d}_i, \bar{a}_{i+1})=1, \ \ i=2,
3, \cdots, s-1.
$$
Also we appoint
$$
\bar{a}_1=\bar{d}_1, \ \ \sum\limits_{i=j}^k(\cdot)=0, \ \ {\rm if}
\ \ k<j
$$
and
$$
\prod\limits_{i=j}^{j-\lambda}(\cdot)=\left\{\begin{array}{l}
1, \ \ \ \ \lambda=1, \\
0, \ \ \ \  \lambda\geq 2.
\end{array}
\right.
$$

\vspace{4mm}

\noindent{\bf Theorem 1.1. (Main Theorem)} \ \ \ {\it If $(a_1, a_2,
\cdots, a_s)|n$, then all integral solutions to the indeterminate
equation (1.1) have the following forms:
$$
\arraycolsep=1.5pt \left\{\begin{array}{rl} x_1=& \displaystyle
n_1\prod\limits_{i=1}^{s-1}\bar{d}_i^{\phi(|\bar{a}_{i+1}|)-1}
+\sum\limits_{m=1}^{s-1}\bar{a}_{m+1}
\prod\limits_{i=1}^{m-1}\bar{d}_i^{\phi(|\bar{a}_{i+1}|)-1}t_m, \\
[3mm] x_k= & \displaystyle
\frac{n_1}{\bar{a}_k}\left(1-\bar{d}_{k-1}^{\phi(|\bar{a}_k|)}\right)
\prod\limits_{i=k}^{s-1}\bar{d}_i^{\phi(|\bar{a}_{i+1}|)-1}
-\bar{d}_{k-1}t_{k-1} \\ [3mm] & \displaystyle
+\sum\limits_{m=2}^{s-1}\frac{\bar{a}_{m+1}}{\bar{a}_k}
\left(1-\bar{d}_{k-1}^{\phi(|\bar{a}_k|)}\right)
\prod\limits_{i=k}^{m-1}\bar{d}_i^{\phi(|\bar{a}_{i+1}|)-1}t_m, \\
[3mm] & \ \ \ \ \ \ \ \ \ \ \  \ \ \ \ \ \ \ \
  k=2, 3, \cdots, s,
\end{array}
\right. \eqno(1.4)
$$
where $t_1, \ t_2, \ \cdots, \ t_{s-1}$ are arbitrary integers. }

\section*{2.  The Proof of Theorem 1.1}

To prove our Theorem 1.1, we restate the following Euler's lemma,
which is required in later analysis.

\bigbreak

\noindent{\bf Lemma 2.1} (Euler's Lemma [2, 3, 5]). \ \ \ {\it Let
$(a,b)=1$. Then
$$
b\left|\left(1-a^{\phi(|b|)}\right)\right.,
\eqno(2.1)
$$
and
$$
a\left|\left(1-b^{\phi(|a|)}\right)\right.,
\eqno(2.2)
$$
where $\phi(\cdot)$ denotes Euler's function.
}

\bigbreak

To study the indeterminate equation (1.1), we first discuss
a simple case of (1.1) with $s=2$.

\bigbreak
\noindent{\bf Lemma 2.2.} \ \ \
{\it
For the indeterminate equation
$$
ax+by=c,
\eqno(2.3)
$$
if $(a,b)|c$, then all integral solutions to (2.3) have the
following forms:
$$
\left\{\begin{array}{l}
\displaystyle
x=c_0a_0^{\phi(|b_0|)-1}+b_0t, \\ [3mm]
\displaystyle
y=\frac{c_0}{b_0}\left(1-a_0^{\phi(|b_0|)}\right)-a_0t
\end{array}
\right.
\eqno(2.4)
$$
or
$$
\left\{\begin{array}{l}
\displaystyle
x=\frac{c_0}{a_0}\left(1-b_0^{\phi(|a_0|)}\right)-b_0t, \\ [3mm]
\displaystyle
y=c_0b_0^{\phi(|a_0|)-1}+a_0t,
\end{array}
\right.
\eqno(2.5)
$$
where $a_0=\frac{a}{(a,b)}$, $b_0=\frac{b}{(a,b)}$,
$c_0=\frac{c}{(a,b)}$, $t=0, \ \pm 1, \
\pm 2, \cdots$.
}

\bigbreak
\noindent{\it Proof.} \ \ \
Without the loss of generality, we only prove (2.4).
The proof of (2.5) is similar and the details are omitted.
In fact, using Lemma 2.1, it is easy to verify that $(x, y)$
is a integral solution to (2.2). On the other hand, let
$(x_0, y_0)$ be a integral solution to (2.2), i.e.,
$$
ax_0+by_0=c.
\eqno(2.6)
$$
Then
$$
a_0x_0+b_0y_0=c_0,
\eqno(2.7)
$$
where $(a_0, b_0)=1$, which implies
$$
a_0x_0\equiv c_0 \ \ \ \ \ ({\rm mod} \ |b_0|). \eqno(2.8)
$$
Therefore, $x\equiv x_0  \ \ ({\rm mod} \ |b_0|)$ must be the
solution of (2.8).

Noticing (2.8) has a unique solution
$$
x\equiv a_0^{\phi(|b_0|)-1}c_0 \ \ \ \ \ ({\rm mod} \ |b_0|),
$$
it follows
$$
x_0\equiv a_0^{\phi(|b_0|)-1}c_0 \ \ \ \ \  ({\rm mod} \ |b_0|).
$$
This shows that there exists a $t_0\in \{0, \pm 1, \pm 2,
\cdots\}$, such that
$$
x_0=c_0a_0^{\phi(|b_0|)-1}+b_0t_0.
\eqno(2.9)
$$
Substituting (2.9) into (2.7), we have
$$
a_0\left(c_0a_0^{\phi(|b_0|)-1}+b_0t_0\right)+b_0y_0=c_0,
\eqno(2.10)
$$
which implies
$$
y_0=\frac{c_0}{b_0}\left(1-a_0^{\phi(|b_0|)}\right) -a_0t_0.
\eqno(2.11)
$$
(2.9) and (2.10) show that every solution $(x_0, y_0)$ to equation (2.3)
satisfies (2.4).

The proof of Lemma 2.2 is completed.

\bigbreak

Now we will seek a formula of all integral solutions to (1.1). To do
this,

\vspace{4mm}

\noindent{\it Proof of Theorem 1.1.} \ \ \ First, we prove that
$(x_1, x_2, \cdots, x_s)$ defined by (1.4) is a integer solution to
(1.1). By using Lemma 1.1, we know that $x_1, \ x_2, \ \cdots, \
x_s$ defined by (1.4) are integers. Moreover, since
$$
n_1a_1=n_1\bar{d}_1d_2=n_1\bar{d}_1\bar{d}_2d_3=\cdots=
n_1\bar{d}_1\bar{d}_2\cdots\bar{d}_{s-1}d_s
=n\bar{d}_1\bar{d}_2\cdots\bar{d}_{s-1},
$$
$$
n_1a_k=n_1\bar{a}_kd_k=n_1\bar{a}_k\bar{d}_kd_{k+1}=\cdots=
n_1\bar{a}_k\bar{d}_k\cdots\bar{d}_{s-1}d_s
=n\bar{a}_k\bar{d}_k\cdots\bar{d}_{s-1},
$$
$$
k=2, 3, \cdots, s-1,
$$
and
$$
n_1a_s=n_1\bar{a}_sd_s=n\bar{a}_s,
$$
we have
$$
\arraycolsep=1.5pt
\begin{array}[b]{rl} & \displaystyle
a_1n_1\prod\limits_{i=1}^{s-1}\bar{d}_i^{\phi(|\bar{a}_{i+1}|)-1}
+a_2\frac{n_1}{\bar{a}_2}\left(1-\bar{d}_1^{\phi(|\bar{a}_2|)}\right)
\prod\limits_{i=2}^{s-1}\bar{d}_i^{\phi(|\bar{a}_{i+1}|)-1}
+\cdots \\ [3mm]
& \displaystyle
+a_{s-1}\frac{n_1}{\bar{a}_{s-1}}\left(1-\bar{d}_{s-2}^
{\phi(|\bar{a}_{s-1}|)}\right)
\prod\limits_{i=s-1}^{s-1}\bar{d}_i^{\phi(|\bar{a}_{i+1}|)-1}
+a_s\frac{n_1}{\bar{a}_s}\left(1-\bar{d}_{s-1}^
{\phi(|\bar{a}_s|)}\right)=n,
\end{array}
\eqno(2.12)
$$
$$
a_1\bar{a}_2t_1-a_2\bar{d}_1t_1=\bar{d}_1d_2\bar{a}_2t_1
-\bar{a}_2d_2\bar{d}_1t_1=0,
\eqno(2.13)
$$
and
$$
\arraycolsep=1.5pt
\begin{array}[b]{rl} & \displaystyle
a_1\bar{a}_{m+1}\prod\limits_{i=1}^{m-1}\bar{d}_i^{\phi(|\bar{a}_{i+1}|)-1}t_m
+a_2\frac{\bar{a}_{m+1}}{\bar{a}_2}\left(1-\bar{a}_1^{\phi(|\bar{a}_2|)}\right)
\prod\limits_{i=2}^{m-1}\bar{d}_i^{\phi(|\bar{a}_{i+1}|)-1}t_m  \\ [3mm]
& \displaystyle
+\cdots
+a_{m-1}\frac{\bar{a}_{m+1}}{\bar{a}_{m-1}}\left(1-\bar{d}_{m-2}^
{\phi(|\bar{a}_{m-1}|)}\right)
\prod\limits_{i=m-1}^{m-1}\bar{d}_i^{\phi(|\bar{a}_{i+1}|)-1}t_m \\ [3mm]
& \displaystyle
+a_m\frac{\bar{a}_{m+1}}{\bar{a}_m}\left(1-\bar{d}_{m-1}^
{\phi(|\bar{a}_m|)}\right)t_m
-a_{m+1}\bar{d}_mt_m  = 0, \\ [3mm]
& \ \ \ \ \ \ \ \ \ \ \ \ \ \ \ \ \ \ \ \ \ \
\ \ \ \ m=2, 3, \cdots, s-1.
\end{array}
\eqno(2.14)
$$
Adding both sides of (2.12), (2.13) and (2.14), we have
$$
a_1x_1+a_2x_2+\cdots+a_sx_s=n,
$$
which implies $(x_1, x_2, \cdots, x_s)$ defined by (1.4) is a
integral solution to (1.1).

On the other hand, we will prove that every integral solution to
(1.1) can be represented into form (1.4) by using induction for $s$.

For $s=2$, it is true by Lemma 2.2.

Suppose that it is true for the indeterminate equation of $s-1$
variables, i.e., the every solution of
$$
a_1x_1+a_2x_2+\cdots+a_{s-1}x_{s-1}=n
$$
can be represented into form (1.4). Now we will show that it is true
for $s$.

Since $d_{s-1}|(a_1x_1+a_2x_2+\cdots+a_{s-1}x_{s-1})$,
there exists $y_{s-1}$ such that
$$
a_1x_1+a_2x_2+\cdots+a_{s-1}x_{s-1}=d_{s-1}y_{s-1}.
\eqno(2.15)
$$
(1.1) and (2.15) show
$$
d_{s-1}y_{s-1}+a_sx_s=n.
\eqno(2.16)
$$
From Lemma 2.2 and the inductive assumption, we have
$$
x_1=
y_{s-1}\prod\limits_{i=1}^{s-2}\bar{d}_i^{\phi(|\bar{a}_{i+1}|)-1}
+\sum\limits_{m=1}^{s-2}\bar{a}_{m+1}
\prod\limits_{i=1}^{m-1}\bar{d}_i^{\phi(|\bar{a}_{i+1}|)-1}t_m,
\eqno(2.17)
$$
$$
\arraycolsep=1.5pt
\begin{array}[b]{rl}
x_k= & \displaystyle
\frac{y_{s-1}}{\bar{a}_k}\left(1-\bar{d}_{k-1}^{\phi(|\bar{a}_k|)}\right)
\prod\limits_{i=k}^{s-2}\bar{d}_i^{\phi(|\bar{a}_{i+1}|)-1}
-\bar{d}_{k-1}t_{k-1} \\ [3mm]
& \displaystyle
+\sum\limits_{m=2}^{s-2}\frac{\bar{a}_{m+1}}{\bar{a}_k}
\left(1-\bar{d}_{k-1}^{\phi(|\bar{a}_k|)}\right)
\prod\limits_{i=k}^{m-1}\bar{d}_i^{\phi(|\bar{a}_{i+1}|)-1}t_m, \\ [3mm]
& \ \ \ \ \ \ \ \ \ \ \  \ \ \ \ \ \ \ \
  k=2, 3, \cdots, s-1,
\end{array}
\eqno(2.18)
$$
$$
y_{s-1}=
n_1\bar{d}_{s-1}^{\phi(|\bar{a}_s|)-1}+\bar{a}_st_{s-1},
\eqno(2.19)
$$
and
$$
x_s=
\frac{n_1}{\bar{a}_s}\left(1-\bar{d}_{s-1}^{\phi(|\bar{a}_s|)}\right)
-\bar{d}_{s-1}t_{s-1}.
\eqno(2.20)
$$
Substituting (2.19) into (2.17), (2.18) and noticing (2.20), we know
every integral solution to (1.1) can be represented into form (1.4).

This completes the proof of Theorem 1.1.

\bigbreak \noindent{\bf Remark 2.4.} \ \ \ {\it The formula (1.4) of
all integral solutions in Theorem 2.3 was deduced from the first
group formula (2.4) of Lemma 2.2. If we use the second group formula
(2.5) to solve the indeterminate equation (1.1), we can obtain the
other formula with different form of all integral solutions to
(1.1). }

\section*{3. Applications}

In this section, we will solve the indeterminate equation (1.3) by
using Theorem 1.1. To do this, we first give Euler's functions
$\phi(2^m)$ and $\phi(3^n)$ as follows:
$$
\left\{\begin{array}{l} \displaystyle
\phi(2^m)=2^m-2^{m-1},  \\ [3mm]
\displaystyle
\phi(3^n)=3^n-3^{n-1}.
\end{array}
\right.
\eqno(3.1)
$$
By applying Theorem 1.1, the all integral solutions to (1.3) can be
represented into the following forms:
$$
\left\{\begin{array}{l} \displaystyle
x=2^{m(3^n-3^{n-1}-1)}+3^nt, \\ [3mm]
\displaystyle
y=\frac{1}{3^n}\left(1-2^{m(3^n-3^{n-1})}\right)-2^mt
\end{array}
\right.
\eqno(3.2)
$$
or
$$
\left\{\begin{array}{l} \displaystyle
x=\frac{1}{2^m}\left(1-3^{n(2^m-2^{m-1})}\right)-3^nt, \\ [3mm]
\displaystyle
y=3^{n(2^m-2^{m-1}-1)}+2^mt,
\end{array}
\right.
\eqno(3.3)
$$
where $t=0, \ \pm 1, \ \pm 2, \cdots$.

\vskip 1cm

{\bf Acknowledgement:} \ \ The research was supported by the Natural
Science Foundation of China $\#$10625105, $\#$11071093, the PhD
specialized grant of the Ministry of Education of China
$\#$20100144110001, and the Special Fund for Basic Scientific
Research  of Central Colleges $\#$CCNU10C01001.

\vskip 1.1cm

\bibliographystyle{plain}

\begin{thebibliography}{99}

\bibitem{A}Apostol, T.M., Introduction to Analytic Number Theory,
Springer-Verlag, Berlin, Heidelberg, New York, 1976.

\bibitem{Dudley}Dudley, U., Elementary Number Theory, W.H. Freeman and Company, New York,
1978.

\bibitem{Gareth}Gareth, A. Jones and J. Mary Jones, Elementary Number Theory,
Springer-Verlag, Berlin, Heidelberg, New York, 1978.

\bibitem{Hardy}Hardy, G.H. and Wright, E.M., An Introduction to the Thoery of
Numbers, Oxford University Press, Walton Street, Oxford OX2 6DP,
1979.

\bibitem{Hua}Hua, L.G., Introduction to Number Theorey, Springer-Verlag,
Berlin, Heidelberg, New York, 1982.

\bibitem{Kenneth}Kenneth Ireland and Michael Rosen, A Classical Introduction to
Modern Number Theory, Springer-Verlag, Berlin, Heidelberg, New York,
1990.

\bibitem{Mathews}Mathews, G.B., Theory of Numbers, New York: Chelsea Publishing,
1961.

\bibitem{Redmond}Redmond, D., Number Theory: An Introduction, New York:
Marcel Dekker, c1996.

\bibitem{Schmidt}Schmidt, W.M., Diophantine Approximations and Diophantine
Equations, Lecture Notes in Math., Vol. 1467, Springer-Verlag,
Berlin, Heidelberg, New York, 1991.



\end{thebibliography}

\end{document}